\documentclass{sig-alternate}
\usepackage{amsfonts,latexsym,amssymb,amsmath,algorithmic,epsfig}
\newlength{\jmr}
\newlength{\barvinok}
\newtheorem{dfn}{Definition}[section]
\newtheorem{rem}[dfn]{Remark}
\newtheorem{prop}[dfn]{Proposition}
\newtheorem{thm}[dfn]{Theorem} 
\newtheorem{lemma}[dfn]{Lemma}
\newtheorem{cor}[dfn]{Corollary}
\newtheorem{ex}[dfn]{Example}

\newtheorem{viro}{Viro's Theorem}

\renewcommand{\mod}{\mathrm{mod}}
\newtheorem{baker}{Nesterenko-Matveev Theorem}

\newtheorem{algor}[dfn]{Algorithm}

\newlength{\smale}
\newcommand{\thth}{^{\text{\underline{th}}}}

\newcommand{\np}{{\mathbf{NP}}}

\newcommand{\feas}{{\mathbf{FEAS}}}
\newcommand{\disc}{{\mathbf{ADISC}_=}}
\newcommand{\sdisc}{{\mathbf{ADISC}_>}}

\newcommand{\gln}{\mathbb{G}\mathbb{L}_n}
\newcommand{\glm}{\mathbb{G}\mathbb{L}_m}

\newcommand{\sat}{\mathbf{3CNFSAT}} 
\newcommand{\pp}{\mathbf{P}}

\newcommand{\nc}{\mathbf{NC}}
\newcommand{\pspa}{\mathbf{PSPACE}}

\newcommand{\eps}{{\varepsilon}}
\newcommand{\Pro}{{\mathbb{P}}}

\newcommand{\barf}{{\bar{f}}}  
\newcommand{\tilf}{{\tilde{f}}}

\newcommand{\R}{\mathbb{R}}

\newcommand{\C}{\mathbb{C}}
\newcommand{\N}{\mathbb{N}}
\newcommand{\Z}{\mathbb{Z}}
\newcommand{\bO}{\mathbf{O}}

\newcommand{\Zn}{\Z^n}

\newcommand{\Rn}{\R^n}

\newcommand{\Cn}{\C^n}

\newcommand{\Cs}{\C^*}
\newcommand{\Rs}{\R^*}
\newcommand{\cA}{{\mathcal{A}}}
\newcommand{\hA}{{\hat{\cA}}}
\newcommand{\cB}{{\mathcal{B}}}
\newcommand{\cC}{{\mathcal{C}}}

\newcommand{\cF}{{\mathcal{F}}}

\newcommand{\cI}{{\mathcal{I}}}

\newcommand{\cT}{{\mathcal{T}}}

\newcommand{\cV}{{\mathcal{V}}}

\newcommand{\Rsn}{{(\R^*)}^n}

\renewcommand{\qed}{$\blacksquare$}
\newcommand{\dia}{$\diamond$}

\newcommand{\newt}{\mathrm{Newt}}

\newcommand{\size}{\mathrm{size}}
\newcommand{\conv}{\mathrm{Conv}}
\newcommand{\aff}{\mathrm{Aff}}

\newcommand{\supp}{\mathrm{Supp}}

\newcommand{\sign}{\mathrm{sign}}

\begin{document}

\conferenceinfo{ISSAC}{'09 Seoul, Korea} 
\title{Faster Real Feasibility via Circuit Discriminants} 
\numberofauthors{3}  
\author{
\alignauthor 
Frederic Bihan\titlenote{Member of the European Research Training Network 
RAAG CT 2001-00271.}\\ 
\affaddr{UFR SFA, Campus Scientifique}\\
\affaddr{73376 Le Bourget-du-Lac Cedex}\\
\affaddr{France}\\
\email{Frederic.Bihan@univ-savoie.fr}\\   
\alignauthor 
J.\ Maurice Rojas\titlenote{
Partially supported by 
NSF individual grant DMS-0211458, NSF CAREER grant 
DMS-0349309, Sandia National Laboratories, and the American 
Institute of Mathematics. }\\ 
\affaddr{TAMU 3368}\\
\affaddr{Department of Mathematics}\\
\affaddr{Texas A\&M University}\\
\affaddr{College Station, Texas \ 77843-3368}\\
\affaddr{USA}\\ 
\email{rojas@math.tamu.edu} 
\alignauthor 
Casey E.\ Stella\titlenote{
Partially supported by NSF grant DMS-0211458.}\\ 
\affaddr{18409 Newell Road}\\
\affaddr{Shaker Heights, OH \ 44122}\\
\affaddr{USA}\\
\email{cestella@gmail.com} 
} 

\date{28 January 2009} 

\maketitle

\noindent
{\scriptsize Rojas dedicates this paper to the memory of his dear 
friend, Richard Adolph Snavely, 1955--2005.} 
\begin{abstract} 
We show that detecting real roots for honestly 
$n$-variate\linebreak $(n+2)$-nomials 
(with integer exponents and coefficients) can be 
done in time polynomial in the {\bf sparse} encoding for any fixed $n$. The 
best previous complexity bounds were exponential in the sparse encoding, even 
for $n$ fixed. 
We then give a characterization of those functions $k(n)$ such that the 
complexity of detecting real roots for $n$-variate \linebreak 
$(n+k(n))$-nomials transitions from $\pp$ to $\np$-hardness as 
$n\longrightarrow \infty$. 
Our proofs follow in large part from a new complexity threshold for 
deciding the vanishing of $\cA$-discriminants of $n$-variate 
$(n+k(n))$-nomials. Diophantine approximation, through linear forms in 
logarithms, also arises as a key tool. 
\end{abstract} 


\section{Introduction and Main Results} 

\noindent
Consider real feasibility: the problem of deciding the existence of real 
roots for systems of polynomial equations. 
In addition to having numerous practical applications 
(see, e.g., \cite{dimacs}), 
real feasibility is an important motivation behind effectivity estimates 
for the Real Nullstellensatz (e.g., \cite{stengle,schmid}), the quantitative 
study of sums of squares \cite{blekherman,swami,reu08}, and their connection to 
semi-definite programming and optimization \cite{parrilo,lassere}. 
In particular, real solving of sparse polynomial systems arises in concrete 
applications such as satellite orbit mechanics \cite{danimarti}, and 
real solving clearly involves real feasibility as an initial step. 
We are thus inspired to derive new algorithms and complexity lower bounds for 
real feasibility, in the refined setting of sparse polynomials. 

To state our results, let us first clarify some basic notation concerning 
sparse polynomials and some well-known 
complexity classes. Recall that $R^*$ is the multiplicative group of 
nonzero elements in any ring $R$. 
\begin{dfn} 
When $a_j\!\in\!\Rn$, the notations $a_j\!=$\linebreak
$(a_{1,j},\ldots,a_{n,j})$, 
$x^{a_j}\!=\!x^{a_{1,j}}_1\cdots x^{a_{n,j}}_n$, and $x\!=\!(x_1,\ldots,x_n)$ 
will be understood. If $f(x)\!:=\!\sum^m_{j=1} c_ix^{a_j}$ 
where $c_j\!\in\!\Rs$ for all $j$, 
and the $a_j$ are pair-wise distinct, then we call $f$ a 
{\bf (real) $\pmb{n}$-variate $\pmb{m}$-nomial}, and we define 
$\supp(f)\!:=\!\{a_1,\ldots,a_m\}$ to be the {\bf support} of $f$. 
We also let $\cF_{n,m}$ denote the set of all $n$-variate 
$\lfloor m\rfloor$-nomials within $\Z[x_1,\ldots,x_n]$. Finally, for any 
$m\!\geq\!n+1$, we let $\cF^*_{n,m}\!\subseteq\!\cF_{n,m}$ 
denote the subset consisting of those $f$ with $\supp(f)$ {\bf not} 
contained in any $(n-1)$-flat. We also call any 
$f\!\in\!\cF^*_{n,m}$ an {\bf honest 
$\pmb{n}$-variate $\pmb{m}$-nomial} (or {\bf honestly $\pmb{n}$-variate}). \dia
\end{dfn}
For example, $1+7x^2_1x_2x^7_3x^3_4-43x^{198}_1x^{99}_2x^{693}_3x^{297}_4$ 
is a $4$-variate trinomial with support contained in a line segment, but it 
has a real root $x\!\in\!\R^4$ iff the honestly univariate trinomial
$1+7y_1-43y^{99}_1$ has a real root $y_1\!\in\!\R$. More generally 
(via Lemma \ref{lemma:hyp} of Section \ref{sub:lin} below), it will be 
natural to restrict to $\cF^*_{n,n+k}$ (with $k\!\geq\!1$) to study the role 
of sparsity in algorithmic complexity over the real numbers. 

We will work with some well-known complexity classes from the 
classical Turing model, briefly reviewed in the Appendix. 
(A more complete introduction can be found in \cite{papa}.)  
In particular, our underlying notion of input size is clarified in Definition
\ref{dfn:size} of Section \ref{sub:cxity} below, and illustrated 
in Example \ref{ex:input}, immediately following our first main theorem.
So for now, let us just recall the basic inclusions 
$\nc^1\!\subseteq\!\pp\!\subseteq\!\np\!\subseteq\!\pspa$.  
While it is known that $\nc^1\!\neq\!\pspa$ the properness of each of the 
remaining inclusions above is a famous open problem. 

\subsection{Sparse Real Feasibility and $\cA$-Discriminant Complexity} 
\begin{dfn}
Let $\R_+$ denote the positive real numbers
and let $\feas_\R$ (resp.\ $\feas_+$) denote the problem of deciding whether an
arbitrary system of equations from $\bigcup_{n\in\N}
\Z[x_1,\ldots,x_n]$
has a real root (resp.\ a root with all coordinates positive).
Also, for any collection $\cF$ of tuples chosen from
$\bigcup_{k,n\in\N}(\Z[x_1,\ldots,x_n])^k$, we let
$\feas_\R(\cF)$ (resp.\ $\feas_+(\cF)$)
denote the natural restriction of $\feas_\R$ (resp.\ $\feas_+$) to inputs in
$\cF$. \dia   
\end{dfn}
\scalebox{.98}[1]{It has been known since the 1980s that 
$\feas_\R\!\in\!\pspa$}\linebreak  
\cite{pspace}, and an $\np$-hardness lower bound was certainly known 
earlier. However, no sharper bounds in terms of sparsity were known earlier in 
the Turing model until our first main theorem. 
\begin{thm}
\label{THM:BIG} 
Let $Z_+(f)$ denote the zero set of $f$ in $\Rn_+$. Then:  
\begin{enumerate} 
\addtocounter{enumi}{-1}
\item{$\feas_+\!\left(\bigcup_{n\in\N} \cF^*_{n,n+1}\right)$ and 
$\feas_\R\!\left(\bigcup_{n\in\N} \cF^*_{n,n+1}\right)$\linebreak 
are in $\nc^1$. 
In particular, when $f\!\in\!\cF^*_{n,n+1}$, $Z_+(f)$ is either empty or 
diffeotopic\footnote{See Definition \ref{dfn:chamber} of Section \ref{sub:disc} 
below.} to $\R^{n-1}_+$, with 
each case actually occuring.} 
\item{For any fixed $n$, $\feas_+(\cF^*_{n,n+2})$ 
and $\feas_\R(\cF^*_{n,n+2})$ are in $\pp$. } 
\item{For any fixed $\eps\!>\!0$, both 
$\feas_+\!\left(\bigcup_{n\in\N}
\cF^*_{n,n+n^\eps}\right)$ and $\feas_\R\!\left(\bigcup_{n\in\N} 
\cF^*_{n,n+n^\eps}\right)$ are $\np$-hard. }  
\end{enumerate} 
\end{thm} 

\noindent
Slightly sharper algorithmic complexity bounds hold when 
we instead work in the BSS model over $\R$ (thus counting arithmetic 
operations instead of bit operations), and this is detailed in 
\cite{prt}. 
\begin{ex}
\label{ex:input}
A very special case of Assertion (1) of Theorem \ref{THM:BIG} implies that 
one can decide --- for any nonzero $c_1,\ldots,c_5\!\in\!\Z$ and 
$D\!\in\!\N$ --- whether \\
\mbox{}\hfill$c_1+c_2x^{999}_1+c_3x^{73}_1x^{19}_3+c_4x^{27D}_2
+c_5x^{74}_1x^{D}_2x_3$
\hfill\mbox{}\\
has a root in $\R^3$, using a number of bit operations 
polynomial in\\
\mbox{}\hfill $\log(D)+\log\left[(|c_1|+1)\cdots (|c_5|+1)\right]$.
\hfill\mbox{}\\ 
The best previous results (e.g., via the critical points method, 
infinitesimals, and rational univariate reduction, as detailed in 
\cite{bpr}) would yield a bound polynomial in\linebreak   
$D+\log\left[(|c_1|+1)\cdots (|c_5|+1)\right]$ 
instead.  \dia
\end{ex}
We thus see that for sparse polynomials, 
large degree can be far less of a complexity bottleneck over $\R$ than
over $\C$. Theorem \ref{THM:BIG} is proved in Section 
\ref{sub:thresh} below.
The underlying techniques include $\cA$-discriminants (a.k.a.\ sparse 
discriminants) (cf.\ Section \ref{sub:disc}), 
Viro's Theorem from toric geometry (see the Appendix, or 
\cite[Thm.\ 5.6]{gkz94}), and effective estimates on linear forms in 
logarithms \cite{baker,nesterenko}. 

In particular, for any collection $\cF_\cA$ of $n$-variate $m$-nomials
with support $\cA$, there is a polynomial $\Delta_\cA$ in the coefficients 
$(c_i)$ called the {\bf $\pmb{\cA}$-discriminant}. Its real zero set 
partitions $\cF_\cA$ into {\bf chambers} (connected components of the
complement) on which the zero set of an
$f\!\in\!\cF_\cA$ has constant topological type.
A toric deformation argument employing Viro's Theorem enables
us to decide whether a given chamber consists of $f$ having empty or 
non-empty $Z_+(f)$. For any $\cA\!\subset\!\Z^n$ of cardinality $n+2$ (in 
sufficiently general position), there is then a compact formula for the 
$\cA$-discriminant that enables us to pick out which chamber contains a given 
$f$: one simply computes the sign of a linear combination of logarithms.
Our resulting algorithms are thus quite implementable, requiring only fast
approximation of logarithms and some basic triangulation combinatorics
for $\supp(f)$.
\begin{ex}
\scalebox{.9}[1]
{Consider $\cA\!:=\!\{(0,0,0),(999,0,0),(73,0,19)$,}\linebreak
$(0,2009,0), (74,293,1)\}$, which gives us
the family of trivariate pentanomials\\
\scalebox{.92}[1]
{$\cF_{\cA}\!:=\!\left.\left\{c_1+c_2x^{999}_1+c_3x^{73}_1x^{19}_3
+c_4x^{2009}_2+c_5x^{74}_1x^{293}_2x_3 \; \right| 
\; c_i\!\in\!\Rs\right\}$.}\linebreak  
Suppose further that $f\!\in\!\cF_\cA$ is an element satisfying
$c_1,c_2,c_3,$\linebreak 
$c_4\!>\!0$ and $c_5\!<\!0$.
It then turns out via Lemma \ref{lemma:ckt} (cf.\ Section \ref{sec:disc} 
below) that $Z_+(f)$ has a degeneracy iff the {\bf $\pmb{\cA}$-discriminant}, 
$\Delta_\cA(c):=$\\ 
\mbox{}\hfill\scalebox{.75}[1]
{$38132829^{38132829}c^{27886408}_1 c^{2677997}_2
c^{2006991}_3 c^{5561433}_4$}\hfill\mbox{}\\
\mbox{}\hspace{.2cm}\scalebox{.75}[1]{$-27886408^{27886408} 2677997^{2677997} 
2006991^{2006991} 5561433^{5561433} c^{38132829}_5$}\\ 
vanishes. In fact, via the techniques underlying Theorem \ref{THM:BIG}, 
$Z_+(f)$ is either empty, a point, or isotopic to a $2$-sphere, 
according as $\Delta_\cA(c)$ is positive, zero,
or negative. Note in particular that determining the sign of
$\Delta_\cA(c)$ is equivalent to determining the sign of\\
\mbox{}\hfill\scalebox{.67}[1]{{\scriptsize
$38132829\log(38132829)+27886408\log(c_1)+2677997\log(c_2)+
2006991\log(c_3)+5561433\log(c_4)$}}\hfill\mbox{}\\
\scalebox{.57}[1]{{\scriptsize $- 27886408\log(27886408) - 
2677997\log(2677997)-2006991\log(2006991)-5561433\log(5561433)
-38132829\log(c_5)$.}} \dia
\end{ex}

While we review $\cA$-discriminants in Section 
\ref{sub:disc} below, it is important to observe now how the 
computational complexity of $\cA$-discriminants closely parallels that of 
$\feas_\R$: compare Theorem \ref{THM:BIG} above with Theorem \ref{THM:DISC} 
below.  
\begin{dfn}
Let $\disc$ (resp.\ $\sdisc$) denote the problem of deciding whether
$\Delta_\cA(f)$ vanishes (resp.\ determining
the sign of $\Delta_\cA(f)$) for an input polynomial  
$f$ with integer coefficients, where $\cA\!=\!\supp(f)$. Finally, 
let $\disc(\cF)$ (resp.\ $\sdisc(\cF)$)
be the natural restriction of $\disc$ (resp.\ $\sdisc$) to inputs 
in some family $\cF$. \dia
\end{dfn}
\begin{thm}
\label{THM:DISC}\mbox{}\\
1.\ $\disc\left(\bigcup\limits_{n\in\N}\cF^*_{n,n+2}
\right)\!\in\!\pp$ and, for any {\bf fixed} $n$,\linebreak  
\mbox{}\hspace{.4cm}$\sdisc(\cF^*_{n,n+2})\!\in\!\pp$.\\ 
2.\ For any fixed $\eps\!>\!0$, both $\disc\!\left(\bigcup\limits_{n\in\N}
\cF^*_{n,n+n^\eps}\right)$ and\linebreak 
\mbox{}\hspace{.4cm}$\sdisc\!\left(\bigcup\limits_{n\in\N}
\cF^*_{n,n+n^\eps}\right)$ are $\np$-hard. 
\end{thm}

\noindent 
Theorem \ref{THM:DISC} is proved in Section \ref{sub:proof2}, after
the development of some necessary theory in Section \ref{sec:back} below.

\subsection{Related Work} 
Earlier work on algorithmic fewnomial theory has \linebreak 
mainly gone in directions 
other than polynomial-time algorithms. For example, Gabrielov and Vorobjov have
given singly exponential time algorithms for weak stratifications
of semi-Pfaffian sets \cite{gv} --- data from which one can  
compute homology groups of real zero sets of a class of functions 
more general than sparse polynomials. Our approach thus highlights a 
subproblem where faster and simpler algorithms are possible. 

Focussing on feasibility, other than the elementary results 
$\feas_\R(\cF_{1,1})\!\in\!\nc^0$ 
and $\feas_\R(\cF_{1,2})\!\in\!\nc^0$, there appear to have been no
earlier complexity upper bounds of the form
$\feas_\R\left(\cF_{1,m}\right)\!\in\!\pp$, or
even $\feas_\R\left(\cF_{1,m}\right)\!\in\!\np$,
for $m\!\geq\!3$. (With the exception of \cite{rojasye}, algorithmic work 
on univariate real polynomials has focussed on algorithms that 
are quasi-linear in the degree. See, e.g., \cite{lickroy}.)   
Echoing the parallels between $\feas_\R$ and 
$\sdisc$ provided by Theorems \ref{THM:BIG} and \ref{THM:DISC}, 
both $\feas_\R\left(\cF_{1,4}\right)\!\stackrel{?}{\in}\!\pp$ and 
$\sdisc(\cF_{1,4})\!\stackrel{?}{\in}\!\pp$ are open problems. 

As for earlier complexity lower bounds for $\feas_\R$ in terms of sparsity,
we are unaware of any. Indeed, it is not even known whether
$\feas_\R(\Z[x_1,\ldots,x_n])$ is $\np$-hard for some fixed $n$. 
Also, complexity lower bounds for the vanishing of discriminants of 
$n$-variate $(n+k(n))$-nomials (with $k$ a slowly growing function of 
$n$) appear to be new. However, recent work shows that the geometry 
of discriminants chambers can be quite intricate already for 
$f\!\in\!\cF^*_{3,3+3}$ \cite{drrs}. 
Also, it was known even earlier that deciding the vanishing of
sparse discriminants of {\bf uni}variate $m$-nomials (with 
$m$ unbounded) is already $\np$-hard
with respect to randomized reductions \cite{ks}. 
Considering Theorems \ref{THM:BIG} and \ref{THM:DISC}, one 
may thus be inclined to conjecture that $\feas_\R(\Z[x_1])$ is 
$\np$-hard. Curiously, over a different family 
of complete fields (the {\bf $\pmb{p}$-adic rationals}), one can
already prove that detecting roots for univariate $m$-nomials 
(with $m$ unbounded) is $\np$-hard with respect to randomized reductions 
\cite{padic1}.

\section{Background and Ancillary\\ Results} 
\label{sec:back} 
After recalling a basic complexity construction, we will present some 
tools for dealing with $n$-variate $(n+1)$-nomials, and then move on to 
$n$-variate $(n+k)$-nomials with $k\!\geq\!2$. {\bf All proofs for the 
results of this section are in the Appendix.} 

\subsection{A Key Reduction}
\label{sub:cxity}
To measure the complexity of our algorithms,  
let us fix the following definitions for input size.
\begin{dfn}
\label{dfn:size}
For any $a\!\in\!\Z$, we define its size, $\size(a)$, to be
$1+\log(1+|a|)$. More generally, we define the {\bf size} of a
matrix $U\!=\![u_{i,j}]\!\in\!\Z^{m\times n}$ to be
$\sum_{i,j} \size(u_{i,j})$.  Also, for any
$f(x)\!=\!\sum^m_{i=1} c_ix^{a_i}\!\in\!\Z[x_1,\ldots,x_n]$,
we define $\size(f)$ to be $\sum^m_{i=1}[\size(c_i)+\size(a_i)]$. 
Finally,
for $F\!=\!(f_1,\ldots,f_k)\!\in\!(\Z[x_1,\ldots,x_n])^k$,
we define $\size(F)\!=\!\sum^k_{i=1}\size(f_i)$. \dia
\end{dfn}

A key construction we will use later in our $\np$-hardness 
proofs is a refinement of an old trick for embedding 
Boolean satisfiability into real/complex satisfiability. 
We refer to the well-known $\sat$ problem, reviewed in 
the Appendix. 
\begin{prop} 
\label{prop:basic} 
Given any $\sat$ instance $B(X)$ with $n$ variables and $N$ clauses, 
let $W_B$ denote\linebreak 
$((\{1\}\times \Pro^1_\C)\cup(\Pro^1_\C\times\{1\}))^{4N-n}$. 
Then there is an \linebreak $(8N-n)\times (8N-n)$ polynomial system $F_B$ with 
the\linebreak 
following properties: 
\begin{enumerate}
\vspace{-.2cm}
\item{$B(X)$ is satisfiable iff 
$F_B$ has a root in \linebreak $\{1,2\}^n\times W_B$.}
\vspace{-.2cm}
\item{$F_B$ has no more than $33N-4n$ monomial terms, \linebreak 
$\size(F_B)\!=\!O(N)$, 
and every root of $F_B$ in\linebreak $(\Pro^1_\C)^{8N-n}$ lies 
in $\{1,2\}^n\times W_B$ and is degenerate.} 
\end{enumerate}

\vspace{-.1cm}
\noindent 
Also, if we define $t_M(z_1,\ldots,z_M)$ to be\linebreak 
$1+z^{M+1}_1+\cdots+z^{M+1}_M 
-(M+1)z_1 \cdots z_M$, then 
\begin{enumerate}
\addtocounter{enumi}{2} 
\item{$t_M$ is nonnegative on $\R^M_+$, with a unique positive root at 
$(1,\ldots,1)$ that happens to be the only degenerate root of $t_M$ in 
$\C^M$.}
\item{If $\eps\!>\!0$, $f\!\in\!\cF^*_{n,n+k}$, and 
$M\!:=\!\left\lceil k^{1/\eps}\right \rceil$, then 
$f(x)+t_M(z)\!\in\!\cF^*_{\eta,\eta+\eta^\delta}$ for 
$\eta\!=\!n+M$ and some positive $\delta\!\leq\!\eps$. In 
particular, $\size(f(x)+t_M(z))\!=\!O\!\left(\size(f)^{1/\eps}\right)$.   
\qed }  
\end{enumerate}
\end{prop} 

\noindent
The seemingly mysterious polynomial $t_M$
defined above will be useful later when we will need to decrease the difference
between the number of terms and variables in certain polynomials.

\subsection{Efficient Linear Algebra on Exponents} 
\label{sub:lin} 
A simple and useful change of variables is to use\linebreak 
monomials in new variables. 
\begin{dfn}
For any ring $R$, let $R^{m\times n}$ denote the set of $m\times n$ matrices
with entries in $R$. For any $M\!=\![m_{ij}]\!\in\!\R^{n\times n}$ 
and $y\!=\!(y_1,\ldots,y_n)$, we define the formal expression 
$y^M\!:=\!(y^{m_{1,1}}_1\cdots y^{m_{n,1}}_n,\ldots, 
y^{m_{1,n}}_1\cdots y^{m_{n,n}}_n)$. We call the substitution 
$x\!:=\!y^M$ a {\bf monomial change of variables}. Also, for 
any $z\!:=\!(z_1,\ldots,z_n)$, we let $xz\!:=\!(x_1z_1,\ldots,x_nz_n)$. 
Finally, let $\gln(\Z)$ denote the group of all 
matrices in $\Z^{n\times n}$ with determinant $\pm 1$ (the set of {\bf 
unimodular} matrices). \dia  
\end{dfn} 

\begin{prop} 
\label{prop:monochange} 
(See, e.g., \cite[Prop.\ 2]{tri}.) 
For any $U,V\!\in\!\R^{n\times n}$, we have the formal identity 
$(xy)^{UV}\!=\!(x^U)^V(y^U)^V$. Also, if 
$\det U\!\neq\!0$, then the function\linebreak 
$e_U(x)\!:=\!x^U$ is an analytic 
automorphism of $\Rn_+$, and preserves smooth points and singular points of 
positive zero sets of analytic functions. Moreover, 
if $\det U\!>\!0$, then $e_U$ in fact induces a diffeotopy on 
any positive zero set of an analytic function. Finally, 
$U\!\in\!\gln(\R)$ 
implies that $e^{-1}_U(\Rn_+)\!=\!\Rn_+$ and that $e_U$ 
maps distinct open orthants of $\Rn$ to distinct open orthants of $\Rn$. \qed 
\end{prop} 

\noindent 
Proposition \ref{prop:monochange}, with minor variations, has been\linebreak  
observed in many earlier works (see, e.g., \cite{tri}). Perhaps the 
only new ingredient is the observation on diffeotopy, which follows easily 
from the fact that $\gln^+(\R)$ (the
set of all $n\times n$ real matrices with positive determinant)
is a connected Lie group. 

Recall that the {\bf affine span}
of a point set $\cA\!\subset\!\Rn$, $\aff \cA$, is the set of real
linear combinations $\sum_{a\in\cA} c_a a$ satisfying
$\sum_{a\in\cA}c_a\!=\!0$.
\begin{lemma}
\label{lemma:hyp}
Given any $f\!\in\!\cF_{n,m}$ with $d=$\linebreak
$\dim \aff(\supp(f))\!<\!n\!<\!m$,
we can find (using a number of bit operations polynomial
in $\size(\supp(f))$) a $U\!\in\!\gln(\Z)$ such that
$g(y)\!:=\!f\!\left(y^U\right)\!\in\!\cF^*_{d,m}$ and $g$
vanishes in $\R^d_+$ (resp.\ $(\Rs)^d$) iff $f$ vanishes in 
$\Rn_+$ (resp.\ $\Rsn$). In particular, there is an
absolute constant $c$ such that $\size(U)\!=$\linebreak
$O(\size(\supp(f))^c)$.
\end{lemma}

To study $Z_+(f)$ when $f\!\in\!
\cF^*_{n,n+1}$ it will help to have a much simpler canonical form.
In what follows, we use $\#$ for set cardinality and 
$e_i$ for the $i\thth$ standard basis 
vector of $\Rn$.  
\begin{lemma}
\label{lemma:can}
For any $f\!\in\!\cF^*_{n,n+1}$ we can compute  
$\ell\!\in\!\{0,\ldots,n\}$ within $\nc^1$ and $\gamma\!\in\!\R+$ such that 
$\barf(x)\!:=\!\gamma+x_1+\cdots+x_\ell-x_{\ell+1}-\cdots -x_n$ 
satisfies: (1) either $f$ or $-f$ has exactly $\ell+1$ positive coefficients, 
and (2) $Z_+\!\left(\barf\right)$ and $Z_+(f)$ are diffeotopic. 
\end{lemma}

\begin{cor} 
\label{cor:signs} 
Suppose $f\!\in\!\cF^*_{n,n+1}$ and 
$\supp(f)\!=$\linebreak
$\{a_1,\ldots,a_{n+1}\}\!\subset\!\Rn$. Then 
\begin{enumerate} 
\item{$f$ has a root in $\Rn_+ \Longleftrightarrow$ not all the coefficients 
of $f$ have the same sign. In particular, $Z_+(f)$ is diffeotopic to 
either $\R^{n-1}_+$ or $\emptyset$. }  
\item{If all the coefficients 
of $f$ have the same sign, then $f$ has a root in 
$\Rsn \Longleftrightarrow$ there are indices $i\!\in\![n]$ and 
$j,j'\!\in\![n+1]$ with $a_{i,j}-a_{i,j'}$ odd. 
} 
\end{enumerate} 
\end{cor} 

\subsection{Combinatorics and Topology of Certain 
$\cA$-Discriminants} 
\label{sub:disc} 
The connection between topology of discriminant complements and
computational complexity dates back to the late 1970s, having been
observed relative to (a) the membership problem for semi-algebraic sets 
\cite{dl79} and (b) the approximation of roots of univariate polynomials 
\cite{smaletopo}. Our goal here is a precise connection between 
$\feas_\R$ and $\cA$-discriminant complements. (See also 
\cite{drrs} for further results in this direction.) 
\begin{dfn} 
\label{dfn:adisc} 
\cite[Ch.\ 1 \& 9--11]{gkz94} 
Given any $\cA\!=\!\{a_1,\ldots,a_m\}\!\subset\!\Zn$ of 
cardinality $m$ and $c_1,\ldots,c_m\!\in\!\Cs$, we 
define $\nabla_\cA\!\subset\!\Pro^{m-1}_\C$ --- 
the {\bf $\cA$-discriminant}\linebreak {\bf variety} --- 
to be the closure of the 
set of all\linebreak $[c_1:\cdots :c_m]\!\in\!\Pro^{m-1}_\C$ such that 
$f(x)\!=\!\sum^m_{i=1} c_ix^{a_i}$ has a degenerate root in $\Cn$. 
We then define $\Delta_\cA\!\in\!\Z[c_1,\ldots,c_m]\!\setminus\!\{0\}$ --- 
the {\bf $\pmb{\cA}$-discriminant} --- to be the unique 
(up to sign) irreducible defining polynomial of $\nabla_\cA$. 
Also, when $\nabla_\cA$ has complex codimension at least $2$, 
we set $\Delta_\cA$ to the constant $1$. 
For convenience, we will sometimes write $\Delta_\cA(f)$ in place of 
$\Delta_\cA(c_1,\ldots,c_m)$. \dia 
\end{dfn} 

To prove our results, it will actually suffice to deal with a small 
subclass of $\cA$-discriminants.
\begin{dfn} 
\label{dfn:ckt} 
We call $\cA\!\subset\!\Rn$ a {\bf (non-degenerate) 
circuit}\footnote{This terminology comes from matroid theory and
has nothing to do with circuits from complexity theory.} 
iff $\cA$ is affinely dependent, but every proper subset of $\cA$ is affinely
independent. Also, we say that $\cA$ is a {\bf degenerate circuit} iff
$\cA$ contains a point $a$ and a proper subset $\cB$ such that $a\!\in\!B$, 
$\cA\setminus a$ is affinely independent, and $\cB$ is a non-degenerate 
circuit. \dia 
\end{dfn}
For instance, both \epsfig{file=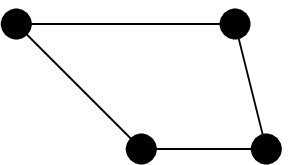,height=.35cm} and
\epsfig{file=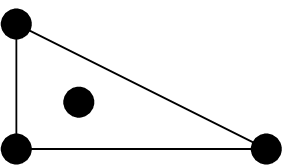,height=.35cm} are circuits, but
\epsfig{file=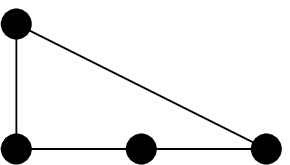,height=.35cm} is a degenerate circuit. 
In general, for any degenerate circuit $\cA$, the subset $\cB$ named above is 
always unique. 

The relevance of $\cA$-discriminants to $m$-nomial zero 
sets can be summarized as follows. 
\begin{dfn} 
\label{dfn:chamber} 
\scalebox{.95}[1]{Following the notation of Definition 
\ref{dfn:adisc},}\linebreak 
we call any connected component of $\Pro^{m-1}_\R\setminus (\nabla_\cA\cup
\{[x_1:\cdots:x_m]\; | \; x_1\cdots x_m\!=\!0\})$ a 
{\bf (real) $\cA$-discriminant chamber}. Also, given any subsets 
$X,Y\!\subseteq\!\Rn_+$, we say that they are {\bf isotopic} (resp.\ 
{\bf diffeotopic}) iff there is a continuous (resp.\ differentiable) function 
$H : [0,1]\times X  
\longrightarrow \Rn_+$ such that $H(t,\cdot)$ is a homeomorphism (resp.\ 
diffeomorphism) for all 
$t\!\in\![0,1]$, $H(0,\cdot)$ is the identity on $X$, and $H(1,X)\!=\!Y$.  
Finally, for any $\cA\!\subset\!\Rn$ of cardinality $m$, let $\cF_\cA$ 
denote the set of all $n$-variate $m$-nomials with support $\cA$. \dia 
\end{dfn} 
\begin{rem} 
Note that when $\cA$ has cardinality $m$, we may naturally identify elements 
of $\Pro^{m-1}_\C$ (resp.\ 
$\Pro^{m-1}_\R$) with equivalence classes determined by 
nonzero complex (resp.\ real) multiples of elements of $\cF_\cA$. \dia 
\end{rem}

The topology of toric real zero sets is known to be constant 
on discriminant chambers (see, e.g., \cite[Ch.\ 11, Sec.\ 5A, Prop.\ 5.2, 
pg.\ 382]{gkz94}). However, we will need a refinement of this fact to 
positive zero sets, so we derive this directly for 
$\cA$ in sufficiently general position --- non-degenerate circuits in 
particular. 
\begin{lemma} 
\label{lemma:chamber} 
Following the notation above, suppose $\cA\!\subset\!\R^n$ is such that 
the minimum of any linear form on $\cA$ is minimized at no more than  
$n+1$ points. Also let $\cC$ be any $\cA$-discriminant chamber. 
Then $f,g\!\in\!\cC \Longrightarrow Z_+(f)$ and $Z_+(g)$ are 
diffeotopic.  
\end{lemma} 

There is then a very compact description for $\nabla_\cA$ when $\cA$ is a 
circuit. 
\begin{lemma}
\label{lemma:ckt} 
Suppose $\cA\!=\!\{a_1,\ldots,a_{n+2}\}\!\subset\!\Zn$ is a non-degenerate 
circuit, $f$ is a polynomial with support $\cA$, $\hA$ is  
the $(n+1)\times (n+2)$ matrix whose $j\thth$ column is $\{1\}\times a_j$, 
$\hA_j$ is the submatrix of $\hA$ obtained by deleting the 
$j\thth$ column, and $b_j\!:=\!(-1)^j\det\hA_j/\beta$ where 
$\beta\!=\!\gcd(\det\hA_1,\ldots,\det\hA_{n+2})$. Then:  
\begin{enumerate}
\item{$\Delta(c_1,\ldots,c_{n+2})$ is, up to a multiple by a nonzero monomial 
term, $\prod\limits^{n+2}_{i=1} \left(\frac{c_i}{b_i}\right)^{b_i}-1$. 
Also, $(b_1,\ldots,b_{n+2})$ can be computed in $\pp$. }  
\item{ $\prod\limits^{n+2}_{i=1} (\sign(b_ic_i)c_i/b_i)^{\sign(b_ic_i)b_i}
\!=\!1$ for some $[c_1:\cdots:c_{n+2}]\!\in\!\Pro^{n+1}_\R$ 
with $\sign(c_1b_1)\!=\cdots=\!\sign(c_{n+2}b_{n+2})   
\Longleftrightarrow Z_+\!\left(\sum^{n+2}_{i=1}c_i x^{a_i}\right)$ 
contains a degenerate point $\zeta$. In particular, 
$Z_+(f)$ has at most one degenerate point. }
\item{$\cA$ has exactly two triangulations: one with 
simplices $\{\conv(\cA\setminus\{b_i\})\; | \; \sign(b_i)\!>\!0\}$, 
and the other with simplices $\{\conv(\cA\setminus\{b_i\})\; | \; 
\sign(b_i)\!<\!0\}$. Moreover, the preceding description also holds 
when $\cA$ is a degenerate circuit.} 
\end{enumerate}
\end{lemma}

\subsection{Complexity of Circuit Discriminants and 
Linear Forms in Logarithms} 
\label{sec:disc} 
Theorem \ref{THM:DISC} is a central tool behind the upper bounds and lower 
bounds of Theorem \ref{THM:BIG}, and is precisely where diophantine 
approximation enters our scenery. To wit, the proof of Assertion (1) of 
Theorem \ref{THM:DISC} makes use of the following powerful result.
\begin{baker}
\cite[Thm.\ 2.1, Pg.\ 55]{nesterenko} 
For any integers $c_1,\alpha_1,\ldots,c_N,\alpha_N$ with
$\alpha_i\!\geq\!2$ for all $i$, define $\Lambda(c,\alpha)\!:=\!
c_1\log(\alpha_1)+\cdots+c_N\log(\alpha_N)$. Then 
$\Lambda(c,\alpha)\!\neq\!0\Longrightarrow 
\log\left|\frac{1}{\Lambda(c,\alpha)}\right|$ is bounded above by\\ 
\mbox{}\hfill {$2.9(N+2)^{9/2}(2e)^{2N+6}
(2+\log\max_j|c_j|)\prod\limits^N_{j=1}\log|\alpha_j|$.\hfill\qed} 
\end{baker}

Assertion (1) of Theorem \ref{THM:DISC} will follow easily from the 
two algorithms we state below, once we prove their correctness 
and verify their efficiency. However, we will first need 
to recall the concept of a gcd-free basis. In essence, a gcd-free 
basis is nearly as powerful as factorization into primes, but is 
far easier to compute.  
\begin{dfn}
\label{dfn:gcd}
\scalebox{.85}[1]{\cite[Sec.\ 8.4]{bs} For any subset 
$\{\alpha_1,\ldots,\alpha_N\}$}\linebreak  
\scalebox{.9}[1]{$\subset\!\N$, a {\bf gcd-free
basis} is a pair of sets $\left(\{\gamma_i\}^\eta_{i=1},
\{e_{ij}\}_{(i,j)\in [N]\times [\eta]}\right)$}\linebreak 
such that
(1) $\gcd(\gamma_i,\gamma_j)\!=\!1$ for all $i\!\neq\!j$, and
(2) $\alpha_i\!=\!\prod^{\eta}_{j=1} \gamma^{e_{ij}}$ for all $i$. \dia
\end{dfn}
\begin{algor}\label{algor:van}\mbox{}\\
{\bf Input:} Integers $\alpha_1,\beta_1,u_1,v_1,\ldots,\alpha_N,\beta_N,
u_N,v_N$. \\
{\bf Output:} A true declaration as to whether
$\alpha^{u_1}_1\cdots\alpha^{u_N}_N\!=\!\beta^{v_1}_1\cdots\beta^{v_N}_N$. \\
{\bf Description:}\\
\vspace{-.5cm}
\begin{enumerate}
\addtocounter{enumi}{-1}
\item{If $\prod^N_{i=1}(\sign \; \alpha_i)^{u_i \ \mod \ 2}\neq
\prod^N_{i=1}(\sign \; \beta_i)^{v_i \ \mod \ 2}$ then
output ``{\tt They are not equal.}'' and {\tt STOP}.}  
\item{Replace the $\alpha_i$ and $\beta_i$ by their absolute values and 
then construct, via Theorem \ref{thm:gcd} of the Appendix, a gcd-free 
basis $(\{\gamma_i\}^\eta_{i=1},\{e_{ij}\}_{(i,j)\in [2N]\times 
[\eta]})$ for $\alpha_1,\ldots,\alpha_N$,
$\beta_1,\ldots,\beta_N$.}
\item{If $\sum^N_{i=1} e_{ij}u_i\!=\!\sum^{2N}_{i=N+1}e_{ij}v_i$ for all
$j\!\in\![\eta]$ then output ``{\tt They are equal.}'' and {\tt STOP}.}
\item{Output ``{\tt They are not equal.}''} 
\end{enumerate}
\end{algor}
\begin{algor}\label{algor:sign}\mbox{}\\
{\bf Input:} Positive integers $\alpha_1,u_1,\ldots,\alpha_M,u_M$ 
and $\beta_1,v_1,$\linebreak
$\ldots,\beta_N,v_N$ with $\alpha_i,\beta_i\!\geq\!2$ for all 
$i$.\\
{\bf Output:} The sign of $\alpha^{u_1}_1\cdots\alpha^{u_M}_M -
\beta^{v_1}_1\cdots\beta^{v_N}_N$. \\
{\bf Description:}\\
\vspace{-.5cm}
\begin{enumerate}
\addtocounter{enumi}{-1}
\item{\scalebox{.9}[1]{Check via Algorithm \ref{algor:van} whether 
$\alpha^{u_1}_1\cdots\alpha^{u_M}_M\!=\!\beta^{v_1}_1\cdots\beta^{v_N}_N$.}
\linebreak 
If so, output ``{\tt They are equal.}'' and {\tt STOP}.} 
\item{Let $U\!:=\!\max\{u_1,\ldots,u_M,v_1,\ldots,v_N\}$, and\\ 
\mbox{}\hfill $E\!:=\!\frac{2.9}{\log 2}(2e)^{2M+2N+6}(1+\log U)$
\hfill\mbox{}\\
\mbox{}\hspace{1.2cm}$\times\left(\prod\limits^M_{i=1}\log|\alpha_i|
\right)\left(\prod\limits^N_{i=1}\log|\beta_i|\right)$. }
\item{For all $i\!\in\![M]$ (resp.\ $i\!\in\![N]$), let $A_i$ (resp.\ $B_i$) 
be a rational number agreeing with $\log\alpha_i$ (resp.\ $\log \beta_i$) in   
its first $2+E+\log_2M$ (resp.\ $2+E+\log_2N$) leading bits.\footnote{For 
definiteness, let us use Arithmetic-Geometric Mean Iteration 
as in \cite{dan} to find these approximations.}}  
\item{Output the sign of
$\left(\sum\limits^M_{i=1}u_iA_i\right)-\left(\sum\limits^N_{i=1}v_iB_i
\right)$ and {\tt STOP}.}
\end{enumerate}
\end{algor}

\begin{lemma}
\label{lemma:bin} 
Algorithms \ref{algor:van} and \ref{algor:sign} are both correct. 
Moreover, following the preceding notation, Algorithms \ref{algor:van} 
and \ref{algor:sign} 
run within a number of bit operations asymptotically linear in, respectively,\\ 
\mbox{}\hfill 
$\sum\limits^N_{i=1} (\log(u_i)\log(\alpha_i)+
\log(v_i)\log(\beta_i))^2$ 
\hfill\mbox{}\\ 
and\\ 
\scalebox{.95}[1]{$(M+N)(30)^{M+N}L(\log U)
\left(\prod\limits^M_{i=1}L(\log(\alpha_i))\right)
\left(\prod\limits^N_{i=1}L(\log(\beta_i))\right)$,}\\
where $L(x)\!:=\!x\log^2(x)\log\log(x)$.  
\end{lemma}

\subsection{Positive Feasibility for Circuits} 
\label{sub:degen} 
For a real polynomial supported on a non-degenerate circuit, there are 
just two ways it can fail to have a positive root: a simple way and a 
subtle way. This is summarized below. Recall that the {\bf Newton polytope of 
$f$} is simply $\newt(f)\!:=\!\conv(\supp(f))$, where $\conv(S)$ denotes 
the {\bf convex hull} (smallest convex set) containing $S$. 
\begin{thm} 
\label{thm:finally}
Suppose $f(x)\!=\!\sum^{n+2}_{i=1}c_ix^{a_i}\!\in\!\cF^*_{n,n+2}$,  
$\supp(f)$ is a non-degenerate circuit, and $b$ is the vector 
from Lemma \ref{lemma:ckt}. Then $Z_+(f)$ is empty 
iff one of the following conditions holds: \\
1.\ All the $c_i$ have the same sign.\\
2.\ $\newt(f)$ is an $n$-simplex and, assuming $a_{j'}$ is the 
unique\linebreak  
\mbox{}\hspace{.4cm}element of $\cA$ lying in the interior of $\newt(f)$,
\linebreak  
\mbox{}\hspace{.4cm}we have $-\sign(c_{j'})\!=\!\sign(c_i)$ for all 
$i\!\neq\!j'$ and\\ 
\mbox{}\hspace{.4cm}$\prod^{n+2}_{i=1} \left(\sign(b_{j'}c_{j'})
\frac{c_i}{b_i}\right)^{\sign(b_{j'}c_{j'})b_i}\!>\!1$.  
\end{thm} 

Positive feasibility for polynomials supported on degenerate circuits can then 
essentially be reduced to the non-degenerate circuit case in some lower 
dimension.  An additional twist arises from the fact that the zero sets of 
polynomials supported on degenerate circuits are, up to a monomial 
change of variables, the graphs of polynomials supported on non-degenerate 
circuits. 
\begin{thm}
\label{thm:degen}
Suppose $f(x)\!=\!\sum^{n+2}_{i=1}c_ix^{a_i}\!\in\!\cF^*_{n,n+2}$ 
has\linebreak 
\scalebox{.94}[1]{support $\cA\!\subset\!\Rn$ that is a degenerate circuit 
with non-degenerate}\linebreak 
subcircuit $\cB\!=\!\{a_1,\ldots,a_{j'}\}$, and 
$b$ is the vector defined in Lemma \ref{lemma:ckt} (ignoring the 
non-degeneracy assumption for $\cA$). Then, when not all the coefficients of 
$f$ have the same sign, $Z_+(f)$ is empty iff both the following conditions  
hold:\\
a.\ $\conv(\cB)$ is a $(j'-2)$-simplex and, permuting indices so\linebreak 
\mbox{}\hspace{.4cm}that $a_{j'}$ is the unique element of $\cB$ lying in the 
relative\linebreak
\mbox{}\hspace{.4cm}interior of $\conv(\cB)$, we have 
$-\sign(c_{j'})\!=\!\sign(c_i)$ for\linebreak 
\mbox{}\hspace{.4cm}all $i\!\neq\!j'$.\\ 
b.\ $\prod^{j'}_{i=1} \left(\sign(b_jc_j)\frac{c_i}{b_i}\right)^{\sign(b_jc_j)
b_i}\!\geq\!1$. 
\end{thm}

\section{The Proofs of Our Main Results: Theorems 
\ref{THM:DISC} and 
\ref{THM:BIG}}  
\label{sec:proofs}

We go in increasing order of proof length. 

\subsection{Proving Theorem \ref{THM:DISC}}
\label{sub:proof2}

\noindent
{\bf Assertion (1):} 
First note
that any input $f$ must have support $\cA\!=\!\{a_1,\ldots,a_{n+2}\}$
equal to either a degenerate circuit or a non-degenerate circuit. Recalling
Assertion (1) of Lemma \ref{lemma:ckt}, observe then that the
vector $b\!:=\!(b_1,\ldots,b_{n+2})$ has a zero coordinate iff
$\cA$ is a degenerate circuit, and $b$ can be computed in time
polynomial in $\size(\cA)$. If $\cA$ is a degenerate circuit then 
(following easily from the definition) $\Delta_{\cA}$ must be identically $1$,  
thus leaving Assertion (1) of our present theorem trivially true.
So let us assume henceforth that $\cA$ is a non-degenerate circuit,
and that $c_j$ is the coefficient of $x^{a_j}$ in $f$ for all $j$.

Via Assertion (1) of Lemma \ref{lemma:ckt} once again, Assertion (1)
of Theorem \ref{THM:DISC} follows routinely from the
complexity bounds from Lemma \ref{lemma:bin}.
In particular, the latter lemma tells us that the bit complexity of
$\disc$, for input coefficients $(c_1,\ldots,c_{n+2})$, is
polynomial in $\sum^{n+2}_{i=1}\log(c_ib_i)$ (following the notation of
Lemma \ref{lemma:ckt}); and the same is true for
$\sdisc$ provided $n$ is fixed.
The classical Hadamard\linebreak  
inequality \cite{mignotte} then tells us that
$\size(b_i)$ is\linebreak $O(n\log(n\max_{j,k}\{a_{jk}\}))$.
So the complexity of $\disc$ is indeed polynomial
in $\size(f)$; and the same holds for $\sdisc$ when $n$ is fixed. \qed

\smallskip
\noindent
{\bf Assertion (2):} We will construct an explicit reduction
of $\sat$ to $\disc\!\left(\bigcup\limits_{n\in\N}
\cF^*_{n,n+n^\eps}\right)$. 
In particular, to any
$\sat$ instance $B(X)$ with $N$ clauses and $n$ variables, let us first 
consider $F_B\!=\!(f_1,\ldots,f_{8N-n})$ --- the associated $(8N-n)\times 
(8N-n)$ polynomial system as detailed in Definition \ref{dfn:sat} of the 
Appendix and Proposition \ref{prop:basic} of Section \ref{sub:cxity}.

Let us then set $M\!:=\!\left\lceil (\max\{0,17N-2n\}
+2)^{1/\eps}\right\rceil$ and define 
the {\bf single} polynomial $f_B$ to be\\
\mbox{}\hfill $f_1+\lambda_1f_2+
\cdots+\lambda_{8N-n-1}f_{8N-n}+\lambda_{8N-n} t_M(z_1,\ldots,z_M)$.\hfill
\mbox{}\\
Letting $\cA$ be the support of $f_B$, it is then easily checked  
(from Definition \ref{dfn:sat} of the Appendix and 
Proposition \ref{prop:basic})  
that $\cA$ is affinely independent and $f_B$ is in 
$\cF^*_{16N-2n+M,N'}$ for some $N'\!\leq\!33N-4n+M+2$. 

By the {\bf Cayley Trick} \cite[Prop.\ 1.7, pp.\ 274]{gkz94} 
we then obtain that $\Delta_\cA(f_B)\!=\!0$ iff\\ 
\mbox{}\hfill ($\star$) $F_B$ has a degenerate root in $(\Pro^1_\C)^{2N-n}$ 
{\bf and}  $t_M$ has\linebreak
\mbox{}\hspace{.6cm}a degenerate root in $(\Cs)^M$.\hfill \mbox{}\\ 
(Since $\newt(t_M)$ is a simplex, it is easily checked that $t_M$ has no 
complex degenerate roots at infinity.) By Proposition \ref{prop:basic}, the 
degenerate roots of 
$F_B$ are exactly $\{1,2\}^n\times W_B$, and $t_M$ has a unique degenerate 
root by construction. So ($\star$) holds iff $B(X)$ has a satisfying 
assignment. We have thus reduced $\sat$ to detecting the vanishing of a 
particular $\cA$-discriminant. 

To conclude, observe that the number of terms of $f_B$ is 
only slightly larger than its number of variables, thanks to 
Proposition \ref{prop:basic}. 
In particular, $\size(f_B)\!=\!O(\size(B)^{1/\eps})$ and 
$f_B\!\in\!\bigcup_{n\in\N}\cF^*_{n,n+n^\delta}
$ for some $\delta\!\in\!(0,\eps]$. Clearly then, 
$\disc\!\left(\bigcup_{n\in\N}\cF^*_{n,n+n^\eps}\right)\!\in\!\pp 
\Longrightarrow \pp\!=\!\np$, thus proving our first desired $\np$-hardness 
lower bound. 

The $\np$-hardness of 
$\sdisc\!\left(\bigcup_{n\in\N}\cF^*_{n,n+n^\delta}\right)$ 
then follows immediately since 
$\sdisc\!\left(\bigcup_{n\in\N}\cF^*_{n,n+n^\delta}\right)$ is a refinement 
of $\disc\!\left(\bigcup_{n\in\N}\cF^*_{n,n+n^\delta}\right)$. \qed 

\subsection{Proving Theorem \ref{THM:BIG}} 
\label{sub:thresh} 

\noindent 
\scalebox{.9}[1]{{\bf Assertion (2):} 
We will give an explicit reduction of $\sat$}\linebreak 
to $\feas_+\!\left(\bigcup_{n\in \N}\cF^*_{n,n+n^\eps}\right)$. 
Attaining such a reduction will require little effort, thanks 
to our earlier reduction used to prove Assertion (2) of Theorem 
\ref{THM:DISC}. 

In particular, for any $\sat$ instance $B$ with $N$ clauses and 
$n$ variables, let us recall the system\linebreak 
$F_B\!=\!(f_1,\ldots,f_{8N-n})$ 
from Definition \ref{dfn:sat} of the Appendix and Proposition \ref{prop:basic}. 
Let us then define $M$ to be\linebreak  
$\left\lceil (\max\{0,42N-n\}+2)^{1/\eps}\right\rceil$ 
and define $g_B(x,z)$ to be\linebreak  
$f^2_1(x)+\cdots+f^2_{4N}(x)+t_M(z_1,\ldots,z_M)$. It is then 
easily checked that $f_B\!\in\!\cF^*_{n+M,N'}$ for some $N'\!\leq\!42N+M+2$. 
Moreover, $B$ has a 
satisfying assignment iff  $g_B$ has a positive root. (Indeed, any root of 
$g_B$ clearly lies in $\{1,2\}^n\times \{1\}^M$.) We have thus reduced 
$\sat$ to a special case of $\feas_+$. 

Now observe that the number of terms of $g_B$ is only 
slightly larger than its number of variables, thanks to Proposition 
\ref{prop:basic}. In particular, $\size(g_B)\!=\!O\!\left(\size(B)^{1/\eps}
\right)$ and $g_B$ is in\linebreak 
$\bigcup_{n\in\N}\cF^*_{n,n+n^\delta}$ for 
some $\delta\!\in\!(0,\eps]$. Clearly then,\linebreak  
$\feas_+\!\left(\bigcup_{n\in \N}\cF^*_{n,n+n^\eps}\right)
\!\in\!\pp \Longrightarrow \pp\!=\!\np$, thus proving one of 
our desired $\np$-hardness lower bounds.  

To conclude, we now need to prove the $\np$-hardness of 
$\feas_\R\!\left(\bigcup_{n\in \N}\cF^*_{n,n+n^\eps}\right)$. 
This we do by employing our preceding argument almost verbatim. 
The only difference is that we instead use the polynomial 
$h_B(x,z)\!:=\!f^2_1(x)+\cdots+f^2_{4N}(x)+ t_M(z^2_1,\ldots,z^2_M)$, 
and observe that $t_M(z^2_1,\ldots,z^2_M)$ is nonnegative 
on all of $\Rn$. So we are done. \qed 

\medskip 
\noindent 
{\bf Assertion (0):} 
Our topological assertion follows immediately from Lemma 
\ref{lemma:can} and Corollary \ref{cor:signs}. 

To obtain our algorithmic assertions, simply note that by 
Assertion (1) of Corollary \ref{cor:signs}, detecting positive roots 
for $f$ reduces to checking whether all the coefficients 
have the same sign. This can clearly be done by $n$ 
sign evaluations and $n-1$ comparisons, doable in 
logarithmic parallel time. So the inclusion involving 
$\feas_+$ is proved. 

Let us now show that we can detect roots in $\Rsn$ within $\nc^1$:  
Employing our algorithm from the last paragraph, we can clearly assume 
the signs of the coefficients of $f$ are all identical (for 
otherwise, we would have detected a root in $\Rn_+$ and finished). 
So then, by Assertion (2) of Corollary \ref{cor:signs}, we can 
simply do a parity check (trivially doable in 
$\nc^1$) of the entries of $[a_2-a_1,\ldots,a_{n+1}-a_1]$. 

To conclude, we simply observe that our algorithm for detecting 
roots in $\Rsn$ trivially extends to root detection in $\Rn$: 
Any root of $f$ in $\Rn$ must lie in some coordinate subspace $L$ of
minimal positive dimension. So, on $L$, the honest $n$-variate 
$(n+1)$-nomial $f$ will restrict to an $f'\!\in\!\cF^*_{n',n'+1}$ with
$n'\!\leq\!n$ and support a subset of the columns of a 
submatrix of $\cA$. So then, we must check whether (a) 
all the coefficients of $f'$ have the same sign or (if not), (b) 
a {\bf sub}matrix of $[a_2-a_1,\ldots,a_{n+1}-a_1]$ has an 
odd entry. In other words, $f$ has a root in $\Rn \Longleftrightarrow 
f$ has a root in $\Rsn\cup\{\bO\}$. Since checking whether $f$ 
vanishes at $\bO$ is the same as checking whether $f$ is missing a 
constant term, checking for roots in $\Rn$ is thus also in $\nc^1$. \qed
\begin{rem}
\label{rem:nc1}
Note that checking whether a given $f\!\in\!\cF_{n,n+1}$ 
lies in $\cF^*_{n,n+1}$ can be done within $\nc^2$: 
one simply finds $d\!=\!\dim \supp(f)$ in $\nc^2$ by 
computing the rank of the matrix whose columns are $a_2-a_1,
\ldots,a_m-a_1$ (via the parallel algorithm of Csanky \cite{csanky}), 
and then checks whether $d\!=\!n$. \dia 
\end{rem}

\smallskip 
\noindent 
\scalebox{.91}[1]{{\bf Assertion (1):} 
The algorithm we use to prove $\feas_+(\cF^*_{n,n+2})$}\linebreak
$\!\in\!\pp$ for fixed $n$ is described just below. Note also that once we 
have $\feas_+(\cF^*_{n,n+2})\!\in\!\pp$ for 
fixed $n$, it easily follows that $\feas_\R(\cF^*_{n,n+2})\!\in\!\pp$: 
The polynomial obtained from an $f\!\in\!\cF^*_{n,n+2}$ by setting any 
non-empty subset of its
variables to $0$ clearly lies in $\cF^*_{n',n'+2}$ for some 
$n'\!<\!n$ (modulo a permutation of variables).
Thus, since we can apply changes of variables like
$x_i\mapsto -x_i$ in $\pp$, and since there are exactly $3^n$
sequences of the form $(\eps_1,\ldots,\eps_n)$ with
$\eps_i\!\in\!\{0,\pm 1\}$ for all $i$, it thus clearly suffices to show
that $\feas_+(\cF^*_{n,n+2})\!\in\!\pp$ for {\bf fixed} $n$.

We thus need only prove correctness, and a suitable complexity bound, 
for the following algorithm: 
\begin{algor}\label{algor:finally}\mbox{}\\
{\bf Input:} A coefficient vector $c\!:=\!(c_1,
\ldots,c_{n+2})$ and a\linebreak (possibly degenerate) circuit 
$\cA\!=\!\{a_1,\ldots,a_{n+2}\}$ of cardinality $n+2$. \\ 
{\bf Output:} A true declaration as to whether $Z_+(f)$ is empty or not, 
where $f(x)\!:=\!\sum^{n+2}_{i=1}c_ix^{a_i}$.\\
{\bf Description:} 
\begin{enumerate}
\item{If all the $c_i$ have the same sign then output 
``{\tt $Z_+(f)\!=\!\emptyset$}'' and {\tt STOP}.} 
\item{Let $b\!=\!(b_1,\ldots,b_{n+2})\!\in\!\Zn$ be the vector obtained by 
applying Lemma \ref{lemma:ckt} to $\cA$. If $b$ or $-b$ has a unique negative 
coordinate $b_{j'}$, and $c_{j'}$ is the unique negative coordinate of $c$ or 
$-c$, then do the following: 
\begin{enumerate}
\item{Replace $b$ by $-\sign(b_{j'})b$, replace $c$ by $-\sign(c_{j'})c$, 
and then reorder $b$, $c$, and $\cA$ by the same\linebreak 
permutation so that $b_{j'}\!<\!0$ and [$b_i\!>\!0$ iff $i\!<\!j'$]. } 
\item{If $j'\!<\!n+2$ and\\ 
\mbox{}\hfill $b^{-b_{j'}}_{j'}\prod^{j-1}_{i=1}c^{b_i}_i\!=\!
c^{-b_{j'}}_{j'}\prod^{n+1}_{i=1}b^{b_i}_i$\hfill\mbox{}\\ 
then output ``{\tt $Z_+(f)\!=\!\emptyset$}'' and {\tt STOP}.} 
\item{Decide via Algorithm \ref{algor:sign} whether\\ 
$b^{-b_{j'}}_{j'}\prod^{{j'}-1}_{i=1}c^{b_i}_i\!\stackrel{?}{>}\!
c^{-b_{j'}}_{j'}\prod^{n+1}_{i=1}b^{b_i}_i$.\\ 
If so, output ``{\tt $Z_+(f)\!=\!\emptyset$}'' and {\tt STOP}.} 
\end{enumerate} }  
\item{Output ``{\tt $Z_+(f)$ is non-empty!}'' and {\tt STOP}.} 
\end{enumerate} 
\end{algor}

\noindent
The correctness of Algorithm \ref{algor:finally} follows 
directly from Theorems \ref{thm:degen} and \ref{thm:finally}. 
In particular, note that $b_i$ is simply the signed volume 
of $\conv(\cA\!\setminus\!\{a_i\})$. So the geometric interpretation 
$b$ or $-b$ having a unique negative coordinate is that the convex hull of 
the unique non-degenerate subcircuit of $\cA$ is a simplex, with  
$a_{j'}$ lying in its relative interior. Similarly, the geometric 
interpretation of $j'\!<\!n+2$ is that $\cA$ is a degenerate 
circuit. Finally, the product comparisons from Steps (b) and (c) simply decide 
the product inequalities stated in  Theorem \ref{thm:finally} 
and Theorem \ref{thm:degen}. 

So now we need only bound complexity, and this follows immediately from Lemma 
\ref{lemma:bin} (assuming we use Algorithm \ref{algor:van} for 
Step (b)). \qed 

It is worth noting that we need to compute the sign of a 
linear combination of logarithms only when the unique non-degenerate 
subcircuit $\cB$ of $\cA$ is a simplex, and all ``vertex'' coefficients have 
sign opposite from the ``internal'' coefficient. Also, 
just as in Remark \ref{rem:nc1}, checking whether a given 
$f\!\in\!\cF_{n,n+2}$ lies in $\cF^*_{n,n+2}$ can be done within $\nc^2$ 
by computing $d\!=\!\dim \supp(f)$ efficiently. Moreover, from our 
preceding proof, we see that deciding whether a circuit is degenerate (and 
extracting $\cB$ from $\cA$ when $\cA$ is degenerate) can be done in $\nc^2$ 
as well, since we can set $\beta\!=\!1$ if we only want the signs of 
$(b_1,\ldots,b_{n+2})$.  

\section*{Acknowledgements}
The authors thank Francisco Santos for earlier discussions on counting regular
triangulations, and Frank Sottile for inviting the second author 
to an April 2008 meeting at the Institute Henri Poincare where a version of 
these results was presented.  
Thanks also to Dima Pasechnik for discussions, and Sue Geller and Bruce 
Reznick for detailed commentary, on earlier versions of this work. 
We also thank AIM and IMA for their hospitality and
support while this paper was nearing completion at respective workshops on 
Random Analytic Surfaces and Complexity, Coding, and Communication. 
Finally, we thank Philippe P\'ebay and David C.\ Thompson for their 
great hospitality at Sandia National Laboratories where this paper was 
completed.  

\bibliographystyle{abbrv}

\section*{Appendix: Complexity, Viro Diagrams, and Postponed Proofs} 

\subsection{Complexity Classes and $\sat$} 
A complete and rigourous description of the complexity
classes we used can be found in \cite{papa}. So for the convenience 
of the reader, we briefly review the following definitions:
\begin{itemize}
\item[$\nc^i$]{ The family of functions computable by Boolean\linebreak 
circuits\footnote{This is the one time we will mention circuits
in the sense of complexity theory: Everywhere else in this paper,
our circuits will be {\bf combinatorial} objects as in Definition
\ref{dfn:ckt}.} with size polynomial in the input size and depth  
$O(\log^i \text{{\tt InputSize}})$. }
\item[$\pp$]{ The family of decision problems that can be done within
time polynomial in the input size. }
\item[$\np$]{ The family of decision problems where a ``{\tt Yes}'' answer can
be {\bf verified} within time polynomial in the input size.}
\item[$\text{\scalebox{.9}[1]{$\pspa$}}$]{ The family of decision problems 
solvable within\linebreak 
polynomial-time, provided a number of processors\linebreak 
exponential in the input size is allowed. }
\end{itemize}

\begin{dfn}
\label{dfn:sat}
Recall that $\sat$ is the problem of deciding whether an $n$-variate
Boolean formula of the form $B(X)=C_1(X)\wedge \cdots \wedge C_N(X)$ has a
satisfying assignment, where each {\bf clause} $C_\ell$ is of one of the
following forms:\\
\mbox{}\hfill $X_i\vee X_j \vee X_k$, \
$X_i\vee X_j \vee \neg X_k$, \hfill\mbox{}\\
\mbox{}\hfill $X_i\vee \neg X_j \vee \neg X_k$,
$\neg X_i\vee \neg X_j \vee \neg X_k$, \hfill \mbox{}\\
$i,j,k\!\in\![n]$ are pairwise distinct,
$\lceil \frac{n}{3}\rceil \!\leq\!N\!\leq\!8\begin{pmatrix} 
n\\ 3\end{pmatrix}$, and a satisfying assigment consists
of an assignment of values from $\{\mathtt{True},\mathtt{False}\}$ to the
variables $X_1,\ldots,X_n$ yielding the equality $B(X)\!=\!\mathtt{True}$
\cite{gj}.
We then define $\size(B)\!:=\!3N$,
$a(x_1,x_2,x_3)\!:=\!(x_1-2)(x_2-2)(x_3-2)$, and $b(x_1)\!:=\!(x_1-1)(x_1-2)$.
Finally, to any $\sat$ clause $C_\ell$ as above, we associate a $4\times 3$
polynomial system $H_{C_\ell}$ as follows: we respectively map clauses of the
form $X_i\vee X_j \vee X_k$, $X_i\vee X_j \vee \neg X_k$, $X_i\vee \neg X_j 
\vee \neg X_k$, $\neg X_i\vee \neg X_j \vee \neg X_k$ to
quadruples of the form\\
\mbox{}\hfill\scalebox{.9}[1]{\begin{tabular}{c}
$(a(x_i,x_j,x_k),b(x_i),b(x_j),b(x_k))$,\\
$(a(x_i,x_j,3-x_k),b(x_i),b(x_j),b(x_k))$,\\
$(a(x_i,3-x_j,3-x_k),b(x_i),b(x_j),b(x_k))$,\\
$(a(3-x_i,3-x_j,3-x_k),b(x_i),b(x_j),b(x_k))$;
\end{tabular}}\hfill\mbox{}\\
and we associate to the $\sat$ instance $B(X)$ a\linebreak 
$(8N-n)\times(8N-n)$ polynomial system with integral coefficients, $F_B$, 
defined to be\\
$(H_{C_1},\ldots,H_{C_N},(u_1-1)(v_1-1),\ldots, 
(u_{4N-n}-1)(v_{4N-n}-1))$. 
In particular, assigning $\mathtt{True}$ (resp.\ $\mathtt{False}$)
to $X_i$ will correspond to setting $x_i\!=\!2$ (resp.\ $x_i\!=\!1$). \dia
\end{dfn}

\noindent
Note that $F_B$ has a natural and well-defined zero set in
$(\Pro^1_\C)^{8N-n}$ since its Newton polytopes are either axes-parallel
line segments or $3$-cubes, and we can multihomogenize with $8N-n$ extra
variables.

\noindent
{\bf Proof of Proposition \ref{prop:basic}:}
Assertions (1) and (2) of Proposition \ref{prop:basic} are elementary. In
particular, the last $4N-n$ polynomials of $F_B$ simply ensure that $F_B$  
has enough variables so that it is square. 
Assertion (3) follows easily from the classical
Arithmetic-Geometric Inequality \cite[Sec.\ 2.5, pp.\ 16--18]{hlp}.
Assertion (4) follows easily upon observing the inequalities
$k\!\leq\!\left\lceil k^{1/\eps}\right\rceil^\eps\!=\!M^\eps\!<\!(n+M)^\eps$
and the fact that $\newt(t_M)$ is $M$-dimensional. \qed

\subsection{Digression on Viro Diagrams} 
\label{sub:viro} 
Let us recall an elegant result of Oleg Viro on the classification 
of certain real algebraic hypersurfaces. In what follows, we liberally 
paraphrase from Proposition 5.2 and Theorem 5.6 of \cite[Ch.\ 5, 
pp.\ 378--393]{gkz94}. \addtocounter{footnote}{1}
\begin{dfn} 
Given any finite point set $\cA\!\subset\!\Rn$, let us call any function 
$\omega : \cA \longrightarrow \R$ a {\bf lifting}, 
denote by $\pi : \R^{n+1}\longrightarrow \Rn$ the natural 
projection which forgets the last coordinate, and let 
$\hat{\cA}\!:=\!\{(a,\omega(a))\; | \; a\!\in\!\cA\}$. We then say that the 
polyhedral subdivision $\Sigma_\omega$ of $\cA$ defined by\\  
$\{\pi(Q) \; | \; Q \text{ a 
lower$^{5}$ facet of } \conv \hat{\cA} \text{ of dimension } 
\dim \aff A\}$ is {\bf induced by the lifting $\omega$}, and we call 
$\Sigma_\omega$ a {\bf triangulation induced by a lifting} iff every cell of 
$\Sigma_\omega$ is a simplex.\footnotetext{A {\bf lower} facet is simply a 
facet which has an inner normal with positive last coordinate.} Finally, 
given any $f(x)\!=\!\sum_{a\in \cA}c_ax^a\!\in\!\Z[x_1,\ldots,x_n]$, we 
define $f_{\omega,\eps}(x)\!:=\!\sum_{a\in \cA}c_a\eps^{\omega(a)}x^a$ to be 
the {\bf toric perturbation of $f$ (corresponding to the lifting $\omega$)}. 
\dia 
\end{dfn} 
\begin{dfn} 
Following the notation above, suppose $\dim \aff \cA\!=\!n$ and $\cA$ 
is equipped with a triangulation $\Sigma$ induced by a lifting {\bf and} a 
function $s : \cA \longrightarrow \{\pm\}$ which we will call a 
{\bf distribution of signs for $\cA$}. We then locally define a 
piece-wise linear manifold  --- the {\bf Viro diagram} $\cV_\cA(\Sigma,s)$ 
--- in the following local manner: For any $n$-cell $C\!\in\!\Sigma$, 
let $L_C$ be the convex hull of the set of midpoints of edges of 
$C$ with vertices of opposite sign, and then define 
$\cV_\cA(\Sigma,s)\!:=\!\bigcup\limits_{C \text{ an } n\text{-cell}} 
L_C$. When $\cA\!=\!\supp(f)$ and $s$ is the corresponding sequence of 
coefficient signs, then we also call $\cV(f)\!:=\!\cV_\cA(\Sigma,s)$ 
the {\bf Viro diagram of $f$}. \dia 
\end{dfn} 
\begin{ex} 
The following figure illustrates $6$ circuits of cardinality $4$, 
each equipped with a triangulation induced by a lifting, and a distribution 
of signs. The corresponding\linebreak 
\scalebox{.97}[1]{(possibly empty) Viro diagrams are drawn 
in thicker lines. \dia}\\ 
\mbox{}\hfill\epsfig{file=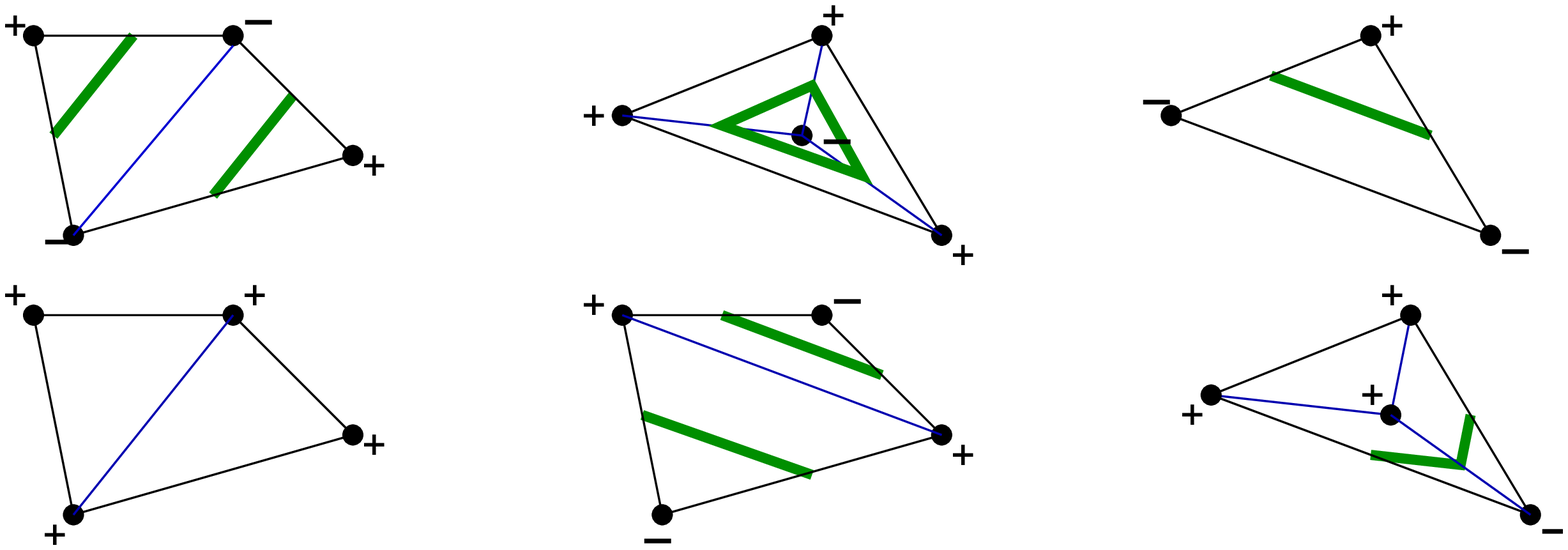,height=1in}\hfill 
\end{ex} 
\begin{viro} 
\label{thm:viro} 
Suppose $f(x)\!=\!\sum_{a\in \cA}c_ax^a$ is in\linebreak 
$\R[x_1,\ldots,x_n]$ with 
$\supp(f)\!=\!\cA$ and $\dim \aff \cA\!=\!n$, $\omega$ is any lifting of 
$\cA$, and define $s_f(a)\!=\!\sign(c_a)$ for all $a\!\in\!\cA$. Then for any 
sufficiently small $\eps\!>\!0$, $Z_+(f_{\omega,\eps})$ is isotopic to 
$\cV_\cA(\Sigma_\omega,s_f)\setminus\partial \conv \cA$. In particular, 
$\cV_\cA(\Sigma_\omega,s_f)$ is a disjoint finite union of piece-wise linear 
manifolds, each possibly having a non-empty boundary. \qed 
\end{viro} 
\begin{lemma} 
\label{lemma:low} 
Suppose $\cA$ is a circuit, $\Sigma$ is a triangulation of $\cA$, 
$n\!=\!\dim \aff \cA$, 
and $s$ is any distribution of signs on $\cA$. 
Then $Z_+(f)$ smooth $\Longrightarrow Z_+(f)$ is isotopic to 
$\cV_\cA(\Sigma,s_f)\setminus\partial \conv \cA$ for some $\Sigma$. 
\end{lemma} 

\noindent 
{\bf Proof of Lemma \ref{lemma:low}:} By Lemma \ref{lemma:ckt} 
it easily follows that $\cA$ 
has at most $2$ discriminant chambers in $\R^{n+2}_+$, and each such 
chamber contains a unique toric perturbation. Since the topology of $Z_+(f)$ 
is constant on any discriminant chamber containing $f$ 
by Lemma \ref{lemma:chamber}, we are done. \qed 
 
Another important consequence of Viro's Theorem deals with how 
``roots at infinity'' sometimes imply the existence of non-compact 
connected components for a positive zero set. 
\begin{lemma} 
\label{lemma:key} 
Assume $f$ has support such that every 
facet of $\newt(f)$ is a simplex and coefficients such that 
$\cV(f)$ intersects $\partial \newt(f)$ for every underlying triangulation. 
Then $Z_+(f)$ has a non-compact connected component. \qed 
\end{lemma} 

\noindent 
Employing our current notation, Lemma \ref{lemma:key} follows directly from 
Lemma 15 of \cite{tri}. 

\subsection{More Postponed Proofs} 

\noindent
{\bf Proof of Lemma \ref{lemma:hyp}:}  
Lemma \ref{lemma:hyp} follows immediately from the following 
well-known factorization for integer matrices and its 
recent complexity bounds. 
\begin{dfn}
\label{dfn:hermite}
\cite{unimod,storjo}
Given any $M\!\in\!\Z^{m\times n}$, the {\bf Hermite factorization}
of $M$ is an identity of the form $UM\!=\!H$ where $U\!\in\!\glm$ and
$H\!=\![h_{ij}]\!\in\!\Z^{n\times n}$ is nonnegative and upper triangular,
with all off-diagonal entries smaller than the positive diagonal entry in
the same column. Finally, the {\bf Smith factorization of $M$} is
an identity of the form $UMV\!=\!S$
where $U\!\in\!\glm$, $V\!\in\!\gln$, and
$S\!=\![s_{ij}]\!\in\!\Z^{m\times n}$ is diagonal, with
$s_{i,i}|s_{i+1,i+1}$ for all $i$. \dia
\end{dfn}
\begin{lemma}
\label{lemma:unimod}
\cite{unimod,storjo}
For any $M\!=\![m_{i,j}]\!\in\!\Z^{n\times n}$, the Hermite and Smith
factorizations of $M$ exist uniquely,
and can be computed within $O(n^4\log^3(n\max_{i,j}|m_{i,j}|))$ bit operations.
Furthermore, in the notation of Definition \ref{dfn:hermite}, the
entries of $U$, $V$, $S$, and $H$ all have bit size\linebreak 
$O(n^3\log^2(2n+\max_{i,j}|m_{i,j}|))$. \qed
\end{lemma}

\medskip 
\noindent
{\bf Proof of Lemma \ref{lemma:can}:} Computing our desired canonical form 
$\barf$ boils down to reordering monomials, performing a monomial change of 
variables, and a rescaling. 

First, let us replace $f$ by $\sign(c_1)f$ and reorder the terms of $f$ 
supported on $\{a_2,\ldots,a_{n+1}\}$ 
so that $c_2,\ldots,c_{\ell'}$ (resp.\ $c_{\ell'+1},\ldots,c_{n+1}$) 
are positive (resp.\ negative), for some unique $\ell'\!\in\!\{1,\ldots,n+1\}$. 
We then form the $n\times n$ matrix $B$ whose 
whose $i\thth$ column is $a_{i+1}-a_1$, for $i\!\in\!\{1,\ldots,n\}$. 
If necessary, let us also  
swap the terms supported on $\{a_2,a_3\}$ (or $\{a_n,a_{n+1}\}$), before 
defining $B$, so that $\det B\!>\!0$ and the sign condition defining $\ell'$ is 
still preserved. (If $n\!=\!1$ then we can simply reorder terms so that 
$a_1\!<\!a_2$ and define $B$ accordingly.) 

Letting $\tilf(x)\!:=\!\frac{f\!\left(x^{B^{-1}}\right)}{x^{B^{-1}a_1}}$ we 
then obtain by Proposition \ref{prop:monochange} that $Z_+(f)$ and 
$Z_+\!\left(\tilf\right)$ are diffeotopic.  Moreover, we clearly have that\\ 
$\tilf(x)\!:=\!c_1+c_2 x_1+\cdots+c_{n+1}x_n$ (remember 
we have permuted the terms of $f$, and thus the $c_i$ as well) where 
$c_1\!>\!0$. So we can now define $\barf(x)$ to be  
$\tilf\!\left(\frac{x_1}{|c_2|},\ldots,\frac{x_n}{|c_{n+1}|}\right)$, 
define $\gamma$ to be the constant term of $\barf$, and 
set $\ell\!:=\!\ell'-1$. It is then easily verified that the coefficients of 
$x_1,\ldots,x_{\ell}$ (resp.\ $x_{\ell+1},\ldots,x_n$) are 
all $1$ (resp.\ $-1$). 

Since $\barf$ was defined by rescaling the variables of $\tilf$, and 
such a scaling of variables can easily be parametrized so as to induce 
a diffeotopy from $Z_+\!\left(\barf\right)$ to $Z_+\!\left(\tilf\right)$,  
Assertion (2) follows. Assertion (1) holds trivially since we never 
altered the difference between the number of positive and negative 
coefficients. That $\ell$ is computable in $\nc^1$ follows from 
Lemma \ref{lemma:unimod}: we can compute $\ell$ simply by sorting, which 
clearly lies in $\nc^1$. \qed 

\medskip 
\noindent
{\bf Proof of Corollary \ref{cor:signs}:} \\
{\bf Assertion (1):}
Employing the canonical form
$\barf(x)\!=\!\gamma+x_1+\cdots+x_\ell-x_{\ell+1}-\cdots-x_n$ of Lemma
\ref{lemma:can} (with $\gamma\!>\!0$ and $\ell\!\in\!\{0,\ldots,n\}$, by
construction), the desired equivalence will follow upon proving that
$\barf$ has a root in $\Rn_+ \Longleftrightarrow \ell\!<\!n$.
The latter equivalence is trivially true. By Lemma \ref{lemma:can}
once more, the statement on diffeotopy type can be reduced to the
special case of $\barf$, which is also immediate. \qed 

\noindent
{\bf Assertion (2):} Dividing by a suitable monomial term,
we can clearly assume that $a_1\!=\!\bO$ and all the coefficients of
$f$ are positive. So it suffices
to prove that $f$ has a root in $\Rsn \Longleftrightarrow$
there are indices $i\!\in\![n]$ and $j\!\in\!\{2,\ldots,n+1\}$ with
$a_{i,j}$ odd. Writing $f(x)\!=\!c_1+c_2x^{a_2}+\cdots+c_{n+1}x^{a_{n+1}}$,
let us now prove the last equivalence.

{\bf ($\Longrightarrow$):} To prove the contrapositive, assume every
$a_{i,j}$ is even.  Then the sign of $f$ is positive on all of $\Rsn$ 
and $f$ thus has {\bf no} roots in $\Rsn$.

{\bf ($\Longleftarrow$):}
Reordering terms and variables, we can clearly assume that $a_{1,1}$ is odd.
Letting $g(x_1)\!:=\!f(x_1,1,\ldots,1)$, note that $g$ must be of the form
$g(x_1)\!=\!c'_1+c'_2x^{a_{1,1}}$, where $c'_1,c'_2\!>\!0$. Since $g$ has
the well-defined real root $-\sqrt[a_{1,1}]{\frac{c'_1}{c'_2}}$,
$f$ then clearly has the root $\left(-\sqrt[a_{1,1}]{\frac{c'_1}{c'_2}},
1,\ldots,1\right)$ which lies in $\Rsn$. \qed

\medskip
\noindent 
{\bf Proof of Lemma \ref{lemma:chamber}:} First, recall that subanalytic sets 
are those sets defined by projections of feasible sets of systems of analytic
inequalities. In particular, $\cC$ is a subanalytic set, and $\cC$ is
thus path connected since $\cC$ admits a decomposition into connected cells.
The existence of such a cell decomposition follows immediately from
the $o$-minimality of subanalytic sets \cite{dries,dries2}. Moreover,
$\cC$ is path connected via {\bf differentiable} paths by the classical
density of $C^\omega$ functions among $C^1$ functions (see, e.g.,
\cite[Ch.\ 2]{hirsch}).

So let $m\!=\!\#\cA$ and let $\phi\!:=\![\phi_1:\cdots:\phi_m] 
: [0,1] \longrightarrow \cC$ be
any differentiable path connecting $f$ and $g$. Also let $\cT$ be 
the positive part of the real toric variety corresponding to $\conv(\cA)$, 
$\cI_\cA\!:=\!\{(c,x)\!\in\!\Pro^{m-1}_\R\times \cT\; | \; 
\sum^m_{i=1} c_ix^{a_i}\!=\!0\}$ the underlying (real) incidence manifold,
and let $\pi$ denote the natural projection mapping
$\Pro^{m-1}_\R\times \cT \longrightarrow \Pro^{m-1}_\R$.
Note then that $\phi$ induces an embedded smooth compact submanifold
$M\!\subseteq\!I_\cA$, consisting of all those $(c,x)$ with
$c\!=\!\phi(t)$ and $\sum^m_{i=1} c_ix^{a_i}\!=\!0$
for some $t$. In particular, we see that $M$ is fibered
over $[0,1]$ and that $\psi\!=\!\phi^{-1}\circ \pi$ is a Morse
function on $M$ with no critical points in $[0,1]$. More to the point,
we obtain a natural flow on $M$ 
inducing a diffeotopy $\delta$ between the zero sets of $f$ and $g$ 
in $\cT$ \cite[Thm.\ 2.2, pg.\ 153]{hirsch}. 

To conclude, we simply observe that the intersection of $\psi^{-1}(t)$ 
with toric infinity is smooth for all $t\!\in\![0,1]$ 
(by our assumption on the facets of $\conv(\cA)$) and thus $\delta$ restricts 
to a diffeotopy between $Z_+(f)$ and $Z_+(g)$. \qed

\medskip 
\noindent
{\bf Proof of Lemma \ref{lemma:ckt}:}
With the exception of the assertion on complexity, Lemma \ref{lemma:ckt}
follows directly from \cite[Prop.\ 1.8, Pg.\ 274]{gkz94},
\cite[Prop.\ 1.2, pg.\ 217]{gkz94}, and the discussion following up to
the end of Section B on page 218 of \cite{gkz94}. In particular, 
the factor $\beta$ takes into account that $\cA$ may not affinely 
generate $\Zn$, but is always the integral affine image of 
an $\cA'$ that is. So the sign condition arises simply from a 
binomial system (with odd determinant) that $\zeta$ must satisfy. 

\scalebox{.95}[1]{The assertion on the complexity of computing $(b_1,\ldots,  
b_{n+2})$}\linebreak 
follows immediately from Csanky's famous parallel algorithm 
for the determinant \cite{csanky}, combined with Lemma \ref{lemma:unimod}.
\linebreak  
Indeed, were it not for the gcd computation for $\beta$, we could instead 
assert an $\nc^2$ complexity bound. \qed 

\medskip 
\noindent 
{\bf Proof of Lemma \ref{lemma:bin}:}  
We first recall the following theorem: 
\begin{thm}
\label{thm:gcd} \cite[Thm.\ 4.8.7, Sec.\ 4.8]{bs} Following the notation of
Definition \ref{dfn:gcd},
there is a gcd-free basis for $\{\alpha_1,\ldots,\alpha_N\}$, with
$\eta$, $\size(\gamma_i)$, and $\size(e_{ij})$
each polynomial in $\sum^N_{\ell=1}\size(\alpha_\ell)$,
for all $i$ and $j$. Moreover, one can always find such a
gcd-free basis using just
$O\!\left(\left(\sum^N_{\ell=1}\size(\alpha_\ell)\right)^2\right)$
bit operations. \qed
\end{thm}

Returning to the proof of Lemma \ref{lemma:bin}, note then that 
Algorithm \ref{algor:van} is correct and runs in the 
time stated by Theorem \ref{thm:gcd} and 
the naive complexity bounds for integer multiplication. 

To prove the remaining half of our lemma, observe 
that the sign of $\alpha^{u_1}_1\cdots\alpha^{u_M}_M-\beta^{v_1}_1\cdots 
\beta^{v_N}_N$ is the same as the sign of 
$S\!:=\!\left(\sum^M_{i=1}u_i\log\alpha_i\right)- 
\left(\sum^N_{i=1}v_i\log\beta_i\right)$.  
Clearly then, $\left|S-\left[\left(\sum^M_{i=1}u_iA_i\right)-
\left(\sum^N_{i=1}v_iB_i\right)\right]\right|\!<\!E/2$ by the 
Nesterenko-Matveev Theorem. So Step (3) of Algorithm \ref{algor:sign} 
indeed computes the sign of $S$ and we thus obtain correctness. 

To see that Algorithm \ref{algor:sign} runs within the time 
stated,  first note that the algorithm computes $M$ (resp.\ $N$) 
approximations of logs of positive integers, each of size 
$O(\max_i \log|\alpha_i|)$ (resp.\ $O(\max_i \log|\beta_i|)$), correct in 
their first $O(E+\log M)$ (resp.\ $O(E+\log N)$) leading bits.  
\noindent
Employing the explicit bit complexity estimates for fast multiplication
of \cite[Table 3.1, pg.\ 43]{bs}, it is easily checked that
Bernstein's method quoted above uses $O(b\log^2(b)\log\log b)$ bit
operations.
So, via our chosen method for approximating logarithms \cite{dan}, 
we see that the complexity of Algorithm \ref{algor:sign} is\\ 
$O\!\left(M(E+\log M)\log^2(E+\log M)\log\log(E+\log M)\right.$\\
\mbox{}\hspace{.4cm}$\left. 
+N(E+\log N)\log^2(E+\log N)\log\log(E+\log N)\right)$\\
$=O\!\left((M+N)E\log^2(E)\log\log E\right)$. 
Upon observing that $M+N\!=\!O\!\left(
\left(\sum^M_{i=1}\log\log|\alpha_i|\right)+ 
\left(\sum^N_{i=1}\log\log|\beta_i|\right)\right)$,\linebreak 
$(2e)^2\!<\!30$, and $abc\log(abc)\!\leq\!a\log(a)b\log(b)c\log(c)$ for 
$a,b,c$ sufficiently large, our final asserted complexity bound follows 
easily. \qed  

\medskip
\noindent
{\bf Proof of Theorem \ref{thm:finally}:} 
First note that Condition (1) implies 
that $f$ maintains the same (non-zero) sign throughout $\Rn_+$. 
So Condition (1) trivially implies that $Z_+(f)\!=\!\emptyset$, and 
we may assume henceforth that not all the coefficients of $f$ have the 
same sign. 
 
Now, if $\newt(f)$ is not a simplex, then every point of $\cA$ 
is a vertex of $\newt(f)$ and thus, independent of the triangulation, any 
Viro diagram for $\cA$ must be non-empty. Since there are only two 
discriminant chambers (by Lemma \ref{lemma:ckt}), Lemma 
\ref{lemma:low} thus implies that $Z_+(f)$ must be non-empty, 
assuming $\Delta_\cA(f)\!\neq\!0$. Lemma \ref{lemma:ckt} 
tells us that $Z_+(f)$ must be non-empty if $\Delta_\cA(f)\!=\!0$. 
So we may assume henceforth that $\newt(f)$ is a simplex and that 
$\Delta_\cA(f)\!\neq\!0$. 

Continuing our focus on Condition (2), note that 
if the sign equalities from Condition (2) fail, then there must exist 
coefficients $c_i$ and $c_{i'}$ of opposite sign such that $a_i$ and 
$a_{i'}$ vertices of $\newt(f)$. So, again, independent of the triangulation, 
any Viro diagram for $\cA$ must be non-empty and thus (just as in the 
preceding paragraph) $Z_+(f)$ must again be non-empty. So we may assume 
henceforth that the sign equalities from Condition (2) hold.  

At this point, it is clear that we need only show that (under our 
current assumptions) $Z_+(f)\!=\!\emptyset 
\Longleftrightarrow$ the discriminant inequality from Condition (2) holds. 
Toward this end, observe that the lifting  
that assigns $a_{j'}\mapsto -1$ and $a_i\mapsto 0$ for all $i\!\neq\!j'$ 
induces the unique triangulation of $\cA$ consisting of a single simplex. 
In particular, the underlying Viro diagram is empty, due to the sign 
equalities. So by Viro's Theorem and Lemma \ref{lemma:low}, $Z_+(f)$ is empty 
for $|c_{j'}|$ sufficiently small. By Lemmata \ref{lemma:chamber} and 
\ref{lemma:ckt}, this topology persists for $|c_{j'}|$ just small enough to 
enforce the discriminant sign stated in Condition (2), so we are done. \qed 

\medskip 
\noindent
{\bf Proof of Theorem \ref{thm:degen}:} 
First, we observe that via an argument almost identical to the
proof of Lemma \ref{lemma:can}, we can find a monomial change of variables
(and multiply by a suitable monomial term) 
so that $\barf(x)\!:=\!x^vf(x^M)$ is of the form\linebreak 
\scalebox{.85}[1]{$c_1+c_2x^{u_1}_1+\cdots
+c_{j'-1}x^{u_{j'-2}}_{j'-2}+c_{j'}x^{\alpha}+c_{j'+1}x^{u_{j'-1}}_{j'-1}
+\cdots+c_{n+2}x^{u_n}_n$,}\linebreak where $u_1,\ldots,u_n\!\in\!\N$ and  
$\alpha\!\in\!\N^{j'-2}\times\{0\}^{n-j'+2}$. 
In particular, defining 
$\barf_\cB(x)\!=\!c_1+c_2x^{u_1}_1+\cdots+c_{j'-1}x^{u_{j'-2}}_{j'-2}
+c_{j'}x^{\alpha}$, it is clear that $Z_+(\barf)$ is nothing more than the 
intersection of the graph of $\barf_\cB$ (which is analytic on $\R^{j'-2}_+$) 
with an orthant. So $Z_+(\barf)$ is smooth and thus, by Proposition 
\ref{prop:monochange}, we obtain that $Z_+(f)$ is diffeotopic to 
$Z_+\!\left(\barf\right)$. It thus suffices to prove our lemma for 
$\barf$.  Note also that the conditions on $Z_+(f)$ allegedly 
characterizing $Z_+(f)\!=\!\emptyset$ are preserved under monomial multiples 
and monomial changes of variables.

Observe now that if $Z_+\!\left(\barf_\cB\right)$ has $\geq\!2$ points then 
$Z_+\!\left(\barf_\cB\right)$
must contain a smooth point $\zeta$, thanks to Lemma \ref{lemma:ckt}. So then,
by the implicit function theorem, $\barf_\cB$ 
must attain both positive and negative values in a neighborhood of 
$\zeta\!\in\!\R^j_+$.
Since the range of $\barf_{B'}\!:=\!c_{j'+1}x^{u_{j'-1}}_{j'-1}+\cdots
+c_{n+2}x^{u_n}_n$ (over $\R^{n-j'+2}$) must contain the positive ray or the 
negative ray, and $\barf\!=\!\barf_\cB-(-\barf_{\cB'})$, we thus obtain that 
$Z_+\!\left(\barf\right)$ is non-empty. 
So we may assume henceforth that $\#Z_+\!\left(\barf_\cB\right)\!\leq\!1$. 

To further simplify matters, observe that if all the coefficients of 
$\barf_\cB$ have the same sign then, by assumption, 
we must have a coefficient of $\barf_{\cB'}$ differing in sign from 
$c_1$. So then, $\barf_\cB(\R^{j'-2}_+)$   
has constant sign and $-\barf_{\cB'}(\R^{n-j'+2})$ contains a 
ray of the same sign.  Therefore $Z_+\!\left(\barf\right)$ is non-empty. We 
may therefore also assume that not all the coefficients of $\barf_\cB$ have 
the same sign. 

We can therefore conclude by studying the two cases\linebreak  
$Z_+\!\left(\barf_\cB\right)\!=\!\emptyset$ and 
$\#Z_+\!\left(\barf_\cB\right)\!=\!1$, combined with all our\linebreak 
assumptions so far. 

\smallskip 
{\bf ($\pmb{Z_+\!\left(\barf_\cB\right)\!=\!\emptyset}$):} 
By Theorem \ref{thm:finally} applied to $\barf_\cB$, we see that Conditions 
(a) and (b) must hold (and the inequality in (b) strictly so), with the 
possible exception 
of the equalities $-\sign(c_{j'})\!=\!\sign(c_i)$ for all $i\!>\!j'$. So we 
need only show that $Z_+(\barf)\!=\!\emptyset \Longleftrightarrow 
-\sign(c_{j'})\!=\!\sign(c_i)$ for all $i\!>\!j'$. 
By an argument almost identical to the last paragraph, we can obtain that 
$\sign(c_{j'})\!=\!\sign(c_i)$ for some 
$i\!>\!j' \Longleftrightarrow \barf_\cB(\R^{j'}_+)$ and 
$-\barf_{\cB'}(\R^{n-j'+2}_+)$ intersect. So the only way 
we can have $Z_+\!\left(\barf\right)\!=\!\emptyset$ 
is for the equalities $-\sign(c_{j'})\!=\!\sign(c_i)$ to hold 
for all $i\!>\!j'$.  

\smallskip 
{\bf ($\pmb{\#Z_+\!\left(\barf_\cB\right)\!=\!1}$):} Clearly, the sole point 
of $Z_+(\barf_\cB)$ must be singular, and thus $\Delta_\cA(f)\!=\!0$, which in 
turn enforces Condition (b) (with equality). 

Now, if $\conv(\cB)$ is not a simplex, then every one of its faces is a 
simplex, since $\cB$ is a non-degenerate circuit. Moreover, all Viro diagrams 
of $\barf_\cB$ (regardless of triangulation) are non-empty and intersect 
$\partial \conv(\cB)$. So, by Lemma \ref{lemma:key}, 
$Z_+(\barf_\cB)$ must contain an unbounded connected component and thus 
$\#\Z_+(\barf_\cB)\!>\!2$. So we may assume that $\conv(\cB)$ is a 
simplex. 
 
At this point, we need only prove that the sign equalities of Condition 
(a) must hold. Toward this end, note that if $\sign(c_{j'})\!=\!\sign(c_i)$  
for some $i\!\in\!\{1,\ldots,j'-1\}$ then all Viro diagrams 
of $\barf_\cB$ (regardless of triangulation) are non-empty and intersect 
$\partial \conv(\cB)$. So, by Lemma \ref{lemma:key} again,  
$Z_+(\barf_\cB)$ must contain an unbounded connected component forcing 
$\#\Z_+(\barf_\cB)\!>\!2$ again. We can therefore assume that 
$-\sign(c_{j'})\!=\!\sign(c_i)$ for all $i\!\in\!\{1,\ldots,j'-1\}$.  

To conclude, we will show that $-\sign(c_{j'})\barf_{\cB}$ attains only 
nonnegative values on $\R^{j'-2}_+$. (This will enforce our final 
desired sign equalities --- $-\sign(c_{j'})\!=\!\sign(c_i)$ for all 
$i\!>\!j'$ --- simply by comparing the range of $\barf_\cB$ and 
$\barf_{\cB'}$ just as before.) First note that  
$\barf_\cB\!\left(\R^{j'-2}_+\right)$ is unaffected by invertible monomial 
changes or positive scalings of variables (thanks to Proposition 
\ref{prop:monochange}). Also, note that $\sign\!\left(\barf_\cB(x)\right)\!=\!
\sign\!\left(x^v\barf_\cB(x)\right)$ for any $x\!\in\!\R^{j'-2}_+$ and 
$v\!\in\!\R^{j'-2}$. 
So it clearly suffices to show that 
$g(x)\!:=\!1+x_1+\cdots+x_{j'-2}-\gamma x^\alpha$\linebreak attains only 
nonnegative values, where now $\alpha\!\in\!\R^{j'-2}_+$ and 
$\alpha_1+\cdots+\alpha_{j'-2}\!<\!1$ 
(since $a_{j'}$ lies in the interior of $\newt(\barf_\cB)$), and 
$\gamma\!=\!\frac{(1-\alpha_1-\cdots-\alpha_{j'-2})^{\alpha_1
+\cdots+\alpha_{j'-2}}}{(1-\alpha_1-\cdots-\alpha_{j'-2})
\alpha^{\alpha_1}_1\cdots \alpha^{\alpha_{j'-2}}_{j'-2}}$ (since Condition 
(b) holds with equality). That $g$ attains only nonnegative values is 
then clearly equivalent to the inequality 
\[ 1+x_1+\cdots+x_{j'-2}\geq \gamma x^\alpha \] 
which is in turn equivalent to\\ 
$(1-\alpha_1-\cdots-\alpha_{j'-2})
(1+x_1+\cdots+x_{j'-2})$\\
\mbox{}\hspace{2.8cm}$\displaystyle{\geq \hspace{.3cm}\prod^{j'-2}_{i=1}
\left(\frac{(1-\alpha_1-\cdots-\alpha_{j'-2})x}{\alpha_i}\right)^{a_i}}$\\ 
or 
\[ (1-\alpha_1-\cdots-\alpha_{j'-2})+
\alpha_1u_1+\cdots+\alpha_{j'-2}u_{j'-2}\geq \prod^{j'-2}_{i=1}
u^{\alpha_i}_i\] 
upon substituting $x_i\!=\!\alpha_iu_i/(1-\alpha_1-\cdots-\alpha_{j'-2})$. 
The last inequality is simply the weighted Arithmetic-Geometric 
Inequality \cite{hlp} so we are done. \qed 
\end{document}